\documentclass[12pt]{article}
\usepackage{amsmath,amssymb}
\usepackage{graphicx}
\newcommand{\bm}[1]{\mbox{\boldmath $#1$}}     

\begin{document}
\newtheorem{prop}{Proposition}
\newtheorem{cor}{Corollary}
\newtheorem{lem}{Lemma}
\newtheorem{thm}{Theorem}

\newcommand{\e}[0]{e}                   
\newcommand{\fb}[0]{\raisebox{.6ex}{\framebox[0.5em]{}}}
\newcommand{\itg}[0]{\mbox{\bf Z}}              
\newcommand{\cpx}[0]{\mbox{\bf C}}              

\title{\Large\bf Quantizing the discrete Painlev\'e VI equation : \\ The Lax formalism}
\author{Koji HASEGAWA
\\
Mathematical Institute, Tohoku University\\
Sendai 980-8578 JAPAN
\thanks{e-mail address : kojihas@math.tohoku.ac.jp
}}
\date{}
\maketitle
\begin{abstract}
A discretization of Painlev\'e VI equation was obtained by Jimbo and Sakai in 1996. 
There are two ways to quantize it:
1) use the affine Weyl group symmetry (of $D_5^{(1)}$) 
\cite{H}, 
\ 2) Lax formalism i.e. monodromy preserving point of view.
It turns out 
 that the second 
approach is also successful and gives the same quantization as in the first approach. 
\end{abstract}

\maketitle
\section{\bf Introduction}\label{sec:intro}

The equation Painlev\'e VI is a 
well known nonlinear ordinary system 
 with rich symmetry and structure.
It can be treated as a  non-autonomous Hamiltonian dynamical system and possesses an extended affine Weyl group symmetry 
of type $D_4^{(1)}$ 
\cite{O}. 
%

The discrete Painlev\'e VI equation ($qP_{VI}$) found by Jimbo and Sakai is the following ordinary difference system: 
we take $t$ as the independent variable (time) of the system and $x(t), y(t)$ the dependent variables.

\begin{equation}{\bf qP_{VI}} \quad 
\begin{cases}
&{\quad
y(t)y(pt)
}
= p^2t^{-2}
\displaystyle
\frac{{x(t)}+a_1^{-2}p^{-1}t}{{x(t)}+a_0^{-2}pt^{-1}}
\cdot
\frac{{x(t)}+a_1^2p^{-1}t}{{x(t)}+a_0^2pt^{-1}},
\\
&{}
\\
\quad
&{x(t)x({p^{-1}}t)}
= t^{-2}
\displaystyle
\frac{{y(t)}+a_4^{-2}t}{{y(t)}+a_5^{-2}t^{-1}}
\cdot
\frac{{y(t)}+a_4^2\;t}{{y(t)}+a_5t^{-1}}.
\end{cases}
\label{eq:qPVI}
\end{equation}
We have five multiplicative parameters, {$p=e^\delta$} : step of time, and  $a_i=e^{\alpha_i}$ ($i=0, 1, 4, 5$). 
The label of the 
parameters are consistently chosen according to the 
${W(D_
5^{(1)})}
$- symmetry :

$s_i(a_j)=a_i^{-C_{ij}}a_j=e^{s_i(\alpha_j)},$
\quad\hfill
$$
[C_{ij}]=\text{the Cartan matrix of type}\ D_5^{(1)}:
{\tiny
\begin{array}{lcccr}
0 &          &   & & 5 \\
  &\backslash&   & \slash   \\
  &          &2-3&          \\
  & \slash   &   & \backslash \\
1 &          &   & &  4 
\end{array}
}
$$

$$
{s_2}(x(t)):=  \frac{a_0a_1^{-1}y(t)+a_2^2}{a_0a_1^{-1}a_2^2y(t)+1} x(t), \quad
{s_j}(x(t)):= x(t) \quad (j\neq 2)
$$

$$
{s_3}(y(t)):=   \frac{a_3^2a_4a_5^{-1}x(t)+1}{a_4a_5^{-1}x(t)+a_3^2} y(t) \quad
{s_j}(y(t)):= y(t) \quad (j\neq 3)
$$
This is the discretized version of the symmetry in the original Painlev\'e VI system ($P_{VI}$) investigated by Okamoto, 
where we adapted the convention used in Tsuda-Masuda \cite{TM}.
The action of the subgroup 
$\langle s_2s_3s_2=s_3s_2s_3, s_i(i=0,1,4,5)\rangle\simeq W(D_4^{(1)})$ 
commutes with the time evolution of ${qP_{VI}}$, and in fact the time evolution itself is a 
translation by the lattice part element $e_3$  which is perpendicular to the 
root lattice $D_4^{(1)}$ embedded in the $D_5^{(1)}$ root lattice (see Appendix). 

In our previous paper \cite{H} we have succeeded in quantizing the affine Weyl group action and thereby construct the 
quantization of the $qP_{VI}$ system. 
Here quantization means the noncommutativity of the dynamical variables 
$x(t), y(t)$ in $qP_{VI}$ and the resulting system $\widehat{qP_{VI}}$ (\ref{qP6^}) looks quite the same as to $qP_{VI}.$ 
Explicit formulae are gathered in the Appendix. 

\medskip
On the other hand, 
one can regard the isomonodromic deformation problem as the origin of the Painlev\'e VI equation. 
The aim of the present paper is to elucidate this point for the quantized discrete equation, that is, 
whether one can obtain $\widehat{qP_{VI}}$
as the quantization of the discretized isomonodoromic deformation problem. 
Actually the answer is quite successful : we obtained the quantization of the Lax form or the Schlesinger equation 
for $\widehat{qP_{VI}}$ (Theorem \ref{thm:C}). 

For this aim, we employed the non-autonomous 
generalization of the quantized lattice system introduced by
Faddeev-Volkov. 
The construction obeys deeply to the quantum group $U_q(A_1^{(1)})$
and its representation; we take the image of the universal R matrix as the Lax matrix or the discrete 
connection matrix.  
The non-autonomous feature comes from the term $c\otimes d$
(where $c$ denotes the canonical central element and $d$ the scaling element, respectively)
 in the universal R matrix of type $A_1^{(1)}$ and 
naturally enters in the pole structure of the Lax matrix, 
which comes from the Heisenberg part of the universal R matirix.


\section{
Review of the Lax formalism for discrete Painlev\'e VI equation}
In this section we review how 
the isomonodoromy deformation problem provides a natural origin of the Painlev\'e VI equation, and how 
one can discretize the problem to obtain qPVI.

Consider the 
 $2\times 2$ regular-
 singular connection on the complex projective line ${\bf P}^1$, 
\begin{equation}
\nabla
=L(z)dz=\displaystyle\sum_{j=1}^n \frac{L^{(j)}}{z-t_j}dz.
\end{equation}
We have $n$ poles $t_1, \cdots, t_n$ at finite points and one at the infinity, put 
$
L^{(\infty)}:=
{\rm Res}_\infty L(z)dz=-
\sum_{j=1}^n L^{(j)}.
$
Let
 $Y(z)$ be the fundamental solution of the linear problem 
$\frac{dY}{dz}=L(z)Y(z)$. 
Then we have the monodoromy matrix $M_j$ along the contour 
$C_j\in \pi_1({\bf P}^1{-}\{t_j\}, *)$   around $t_j$, 
where $*$ stands for the fixed base point (which is different from the singularity): 
$$
{C_j}_*(Y)(z){=}Y(z){M_j}.
$$ 
The matrix 
$M_j$ is conjugate to $e^{2\pi iL^{(j)}}$ and 
satisfy the relation 
$M_1\cdots M_n M_\infty=1$. 

\bigskip
{\bf Fact. } 
The monodromy matrices 
$\{M_j\}$ are constant ({\sl isomonodromy}) with respect to $t_j$'s 
 if the following relations 
hold :
\begin{equation}
\frac{\partial Y}{\partial t_j}Y^{-1}
=-\frac{L^{(j)}}{z-t_j}(=:B_j) \qquad (j=1, \cdots, n).
\label{eq:isomonodromy}
\end{equation}
If this is the case, the compatibility of 
(\ref{eq:isomonodromy}), called the {\sl Schlesinger equation}, should be satisfied: 
\begin{equation}
{[\frac{\partial}{\partial z}-L(z,\vec{t}), \frac{\partial}{\partial t_j}-B_j(z,\vec{t})]=0}
\quad (i, j=1, \cdots, n), 
\label{eq:Schlesinger}
\end{equation}
where $\vec t:=(t_1,\cdots, t_n)$ and $\vec t$ dependence of $L$ and $B_j$ are explicitly written. 
See Jimbo-Miwa-Ueno \cite{JMU} for details. 

This is the Lax form of the isomonodromy problem. 
The case $n+1=4$ reproduces the Painlev\'e VI equation: 
one can assume 
$
(t_1, t_2, t_3, \infty)=(0,1,{t},\infty)
$
and take the dependent variable 
$y(t)$ to be (roughly speaking) the off-diagonal 
 element of $L^{(3)}$. 








According to Jimbo and Sakai \cite{JS}, 
the difference equation case goes quite similarly. 
Let us consider the difference equation 
$$
\frac{Y(qz)-Y(z)}{qz-z}=L(z)Y(z)
\quad
(L(z)=\frac{L^{(1)}}{z}+\frac{L^{(2)}}{z-1}+\frac{L^{(3)}}{z-t}, \ \text{generic})
$$
which can be rewritten as 
$
Y(qz)=\{1+(q-1)zL(z)\}Y(z)
$. 
There exists some function $\gamma$ 
such that 
$(z-1)(z-t)=
{\gamma(qz)}{\gamma(z)}^{-1}$. 
Put ${\cal Y}=\gamma Y$, then
we have
${{\cal Y}(qz)={\cal L}(z){\cal Y}(z)}$, where 
$
{\cal L}(z):=(z-1)(z-t)\{1+(q-1)zL(z)\}
$
is polynomial in $z$. 
Now singularities are $0$ and $\infty$; $1$ and $t$ can be detected as the zero of $\det {\cal L}$.

There exists an solution at $z=0$ of the form 
${{\cal Y}_0(z)}=z^{{\cal L}(0)}\times$(power series in $z$), 
and similarly, ${{\cal Y}_{\infty}(z)}$ for $z=\infty$. 
The `connection' matrix 
\begin{equation}
{{\cal M}(z)}:={\cal Y}_0(z)^{-1}{\cal Y}_\infty(z), \quad {\cal M}(z)={\cal M}(qz), 
\label{eq:connectionmatrix}
\end{equation}
plays the role of the monodoromy matrix. 
Deformation preserving condition 
${\cal M}(z)={\cal M}(z,t)$ 
is then satisfied if we have some ${\cal B}(z, t)$ such that 
${\cal Y}(z, qt)={\cal B}(z,t){\cal Y}(z,t)$. 

The compatibility now reads as {\sl discrete Schlesinger equation},
\begin{equation}
{\cal L}(z,qt){\cal B}(z,t)={\cal B}(qz,t){\cal L}(z,t) 
\label{eq:qSchlesinger}
\end{equation} 
from which Jimbo and Sakai derived the $qP_{VI}$ equation (\ref{eq:qPVI}). 

\medskip
Our goal will be the quantization of (\ref{eq:qSchlesinger}) 
as well as to confirm that it reproduces the quantization 
$\widehat{qP_{VI}}$ (\ref{qP6^}) of ${qP_{VI}}$ (\ref{eq:qPVI}).


\section{The quantized local Lax matrix}


For our aim, we use 
non-autonomous modification of 
Faddeev-Volkov quantization of discrete sine-Gordon equation and its periodic reduction.
In this section we will give the local Lax matrix, which 
can be said as nonautonomously modified 
Izergin-Korepin Lax matrix \cite{IK}.

\bigskip
Let $q$ be a complex number with $0<|q|<1$. 
For $\pm=+$ or $-$ respectively, 
let 
$U_q^\pm=U_q^\pm(A_1^{(1)})$ be the upper/lower subalgebra of the quantum group 
$U_q=U_q(A_1^{(1)})$ generated by the upper/lower Chevalley generators $e^\pm_i$ 
together with the Cartan part $h_i (i=0, 1), d$, where $d$ is the scaling element.
We write the canonical central element as $c(:=h_0+h_1)$. 

Let $c^\pm\in {\bf C}$, and let 
$\rho^\pm$ 
be the representation of $U_q^\pm$ on the space
$V^\pm:=
{\bm C}[e^{\pm\alpha_0}, e^{\pm\alpha_1}]
$ 
defined by
\begin{equation}
e_i^\pm\mapsto -(q-q^{-1})e^{\pm \alpha_i}=:E_i^\pm, \ h_i\mapsto h_i\ (i=0, 1)
, \ c\mapsto c^\pm\in{\bm C}
\label{def:rhopm}
\end{equation}
respectively. 
By the definition $h_i$ acts as the derivation satisfying 
$
[h_i, e^{\pm\alpha_j}]=\pm\alpha_j(h_i)e^{\pm\alpha_j}.
$
%

Let  ${\cal R}\in U_q^+\otimes U_q^-$ be the universal $R$ matrix of $U_q$ and 
$\raisebox{4pt}{{\fbox{}}}_{z}$ 
be the two dimensional evaluation representation of $U$. 
Write 
$k:=q^{h_1}, \Delta^\pm=E_0^\pm E_1^\pm$.  
We have :
$$
kE^{\pm}_1k^{-1}=q^{\pm 2}E_1,\ 
kE^{\pm}_0k^{-1}=q^{\mp 2}E^\pm_0, \
[E^+_1, E^-_1]=(k-k^{-1})(q-q^{-1})
$$ 
and also
\begin{equation}
q^d
\Delta^{\pm}q^{-d}{=}q^{\pm1}\Delta^{\pm}.
\label{eq:Deltamoves}
\end{equation}
Other than (\ref{eq:Deltamoves}), $\Delta^\pm$ commutes with the generators $e^\pm_i, h_i \ (i=0, 1)$ 
of $U_q^\pm$, i.e. 
$
\Delta^\pm\in {\cal Z}({U_q^\pm}'), 
$ 
the center of the derived algebra of $U_q^\pm$.
Put
\begin{equation}
{L_z^+(
\Delta^+)}
{:=}
(\rho^+\otimes\raisebox{0.5ex}%
{{\fbox{}}}_{z})({\cal R}), 
\quad
{L_z^-(\Delta^-)}{:=}
(\raisebox{1.5pt}{{\fbox{}}}_{z}\otimes\rho^-)({\cal R}).
\label{def:Lpm}
\end{equation}
These are the local Lax matrices for our aim. 

\begin{prop} We have
\begin{align}
{L_z^+(
\Delta^+)}
&=
\frac{(q^4z^{-1}\Delta^+,q^4)_\infty}{(q^2z^{-1}\Delta^+,q^4)_\infty}
{
\left[\begin{array}{cc}
 1   &  \frac{1}{z} E^+_0\\
 E^+_1 &   1
\end{array}\right]
\left[\begin{array}{cc}
k^{-\frac{1}{2}} & 0 \\ 0 & k^\frac{1}{2}
\end{array}\right]
}
q^{-c^+d}, 
\label{eq:L^+}
\\
{L_z^-(\Delta^-)}
&=
\frac{(q^4z\Delta^-,q^4)_\infty}{(q^2z\Delta^-,q^4)_\infty}
{
\left[\begin{array}{cc}
  1  &   E^-_1 \\zE^-_0  &   1
\end{array}\right]
\left[\begin{array}{cc}
k^{-\frac{1}{2}} & 0 \\ 0 & k^\frac{1}{2}
\end{array}\right]
}
q^{-c^-d}, 
\label{eq:L^-}
\end{align}
where we used the standard notation for the infinite product : 
$(x,Q)_\infty:=\prod_{n=0}^\infty (1-xQ^n).$
\label{prop:formulaforL}
\end{prop}
This is derived from the formula for the universal R matrix (see e.g. 
\cite{DKP}) 
and the above definition of the corresponding representations. 
The case $c^\pm=0$ is essentially the one used in \cite{IK} and later reproduced by \cite{FV}. 
The infinite product factors come from the Heisenberg part 
(contribution from the null root vectors) in the product formula of ${\cal R},$ 
and therefore we have infinite poles here. 
The pole location will move 
according to (\ref{eq:Deltamoves})  during the dynamics defined in the next section ; 
this is the non-autonomous nature of our 
Lax matrix. 
In the $q\rightarrow 1$ limit, we have 
\begin{equation}
\lim_{q\rightarrow 1}
\left(\frac
{(q^4z^{-1}\Delta^+,q^4)_\infty
}{(q^2z^{-1}\Delta^+,q^4)_\infty
}\right)^{-2}
=
1-\frac{\Delta^+}{z}
\label{eq:limit_of_infiniteproduct}
\end{equation}
%
which is equal to the determinant of the matrix 
\begin{equation}
\left[\begin{array}{cc}
 1   &  \frac{1}{z} E^+_0\\
 E^+_1 &   1
\end{array}\right].
\end{equation}
Similar relation holds for $L^-$, 
showing the $SL(2)$ nature of these local Lax matrices.

Note also 
\begin{equation}
q^d L^\pm_z(\Delta^\pm)q^{-d}
=L^\pm_{zq^\mp 1}(\Delta^\pm), 
\label{eq:homogeneity}
\end{equation}
which can be easily seen from (\ref{eq:L^+}) and (\ref{eq:L^-}).

\bigskip
Let
$
R(\Delta^+, \Delta^-):=(\rho^+\otimes\rho^-)({\cal R}). 
$
According to the multiplicative formula of the universal R matrix, 
it is explicitly written in terms of the quantum dilogarithm \cite{Kir}, 
\begin{equation}
R(\Delta^+, \Delta^-)
=(qE^+_1\otimes E^-_1, q^2)_\infty^{-1}
(q^2\Delta^+\otimes\Delta^-,q^4)_\infty^{-1}
(qE^+_0\otimes E^-_0, q^2)_\infty^{-1}
q^{-T}, \ 
\end{equation}
where
$
T:=\frac{1}{2}h_1\otimes h_1+c^+\otimes d 
+ d \otimes c^-.
$

The Yang-Baxter equation for the universal R matrix
$$\mathcal R_{12}\mathcal R_{13}\mathcal R_{23}
=
\mathcal R_{23}\mathcal R_{13}\mathcal R_{12}
\in U_q^{+}\otimes U_q\otimes U_q^-
$$
immediately implies the following (although it can be checked directly):
\begin{prop}
\begin{equation}
{
L_z^+(\Delta^+)}{R(\Delta^+,\Delta^-)}{L_z^-(\Delta^-)}
=
{L_z^-(\Delta^-)}{R(\Delta^+,\Delta^-)}{L_z^+(\Delta^+)}.
\label{prop:theYBE}
\end{equation}
\label{prop:LRL=LRL}
\end{prop}
\begin{center}
{\sl Figure 1. The Yang-Baxter equation (\ref{prop:theYBE}).}
\\
\bigskip
\includegraphics[height=3cm]{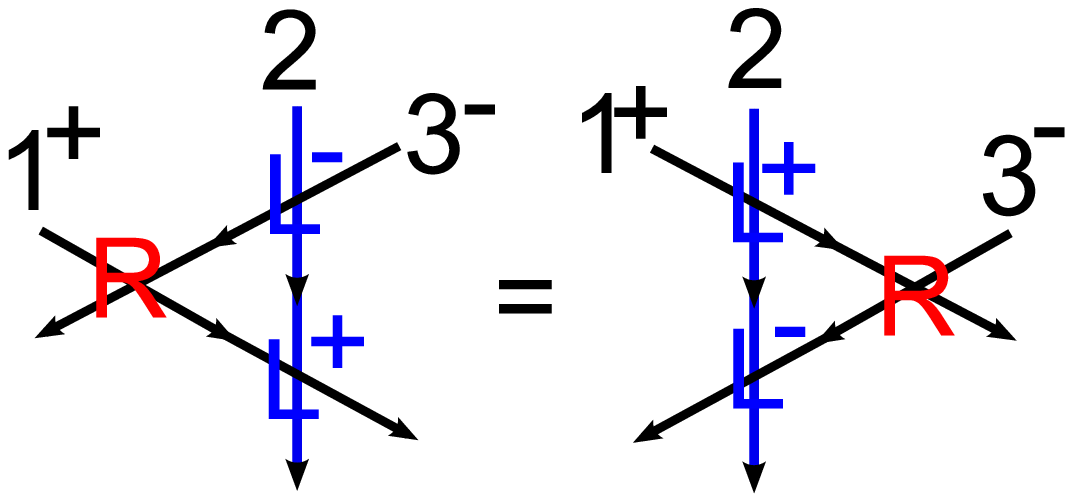}
\end{center}

We can rewrite the above Yang-Baxter relation (Figure 1) among $L^\pm$ and $R(\Delta^+, \Delta^-)$ as follows 
({\sl exchange dynamics} of $L^+$ and $L^-$, Figure 2):
\begin{center}
{\sl Figure 2. The exchange dynamics (\ref{eq:localexchangedynamics}).}
\\
\bigskip
\includegraphics[height=4cm]{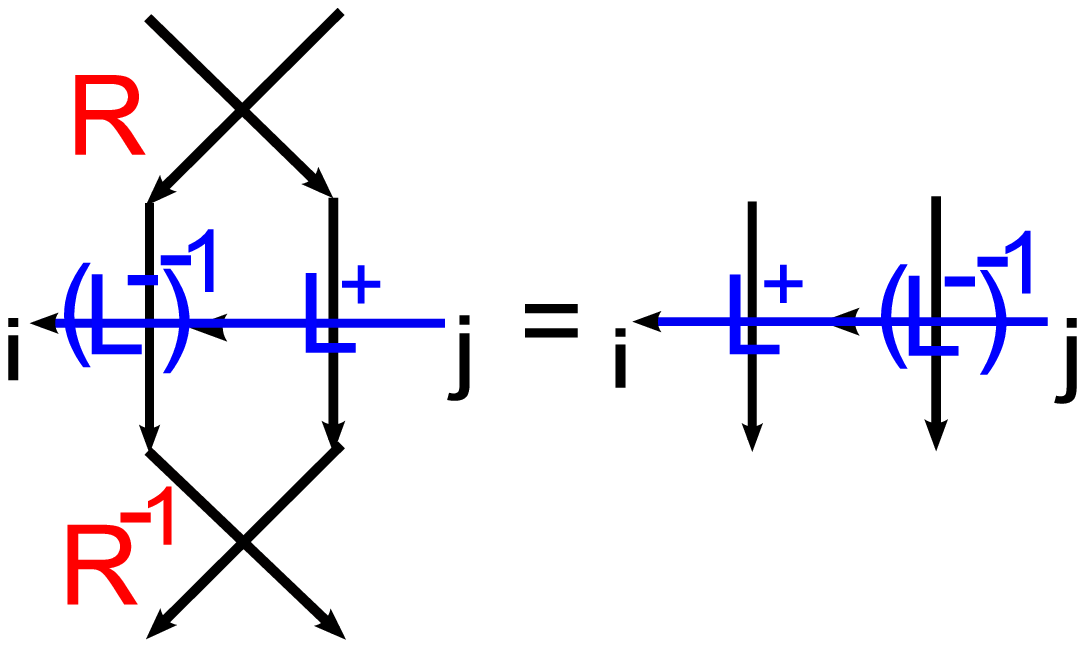}
\end{center}
\begin{equation}
{R^{-1}}\left[{(L_z^-)^{-1}L_z^+}\right]_{ij}{R}
 =\left[{L_z^+(L_z^-)^{-1}}\right]_{ij} 
\quad \in (\rho^+\otimes\rho^-)(U_q^+\otimes U_q^-)
\label{eq:localexchangedynamics}
\end{equation}
where  $[\cdot ]_{ij}$ stands for the matrix element. 
If we work with $U_q(A_\ell^{(1)})$, one can see that this gives the 
quantized periodic Toda lattice equation (cf. \cite{jimbo}). 

\section{The discrete time lattice dynamics and its reduction}

Let us define the one dimensional lattice system. 
Let us write the even/odd lattice points  of $\bm Z$ as $n^+, n^-$, rather than $2n, 2n+1$. 
Then we attach the representation
$$
\rho^\pm=\rho^{n\pm} : U_q^\pm\rightarrow {\rm End}(V^{n\pm}), 
\ V^{n\pm}\simeq {\bf C}[e^{\pm\alpha_0}, e^{\pm\alpha_1}]
$$ 
specified by the parameters 
 $c^{n\pm}:=\rho^{n\pm}(c), \Delta^{n\pm}:=\rho^{n\pm}(\Delta^\pm)$,  
for each of these points $n^\pm$. 
Let 
\begin{equation}
{\cal V}:=\otimes_{n, \pm} V^{n\pm}, \qquad
{\cal A}:=\otimes'_{n, \pm} {\rm End}(V^{n\pm}), 
\end{equation}
where $\otimes'$ stands for the restricted tensor product with respect to $1$. 
${\cal A}$ is the algebra finitely generated by the local operators on $V^{n\pm}$
(i.e. elements of ${\rm End}(V^{n\pm})$). 

We write the $R$ matrix $\rho^{m+}\otimes\rho^{n-}({\cal R})$ as $R^{m+, n-}$ for short and
define the Hamiltonian of our dynamics to be 
${\cal H}={\cal H}_0{\cal H}_1$, 
where (cf. \cite{FV})
$$
{\cal H}_0{:=}\cdots{R^{1^+1^-}R^{2^+2^-}}\cdots, \quad
{\cal H}_1{:=}\cdots{R^{2^+1^-}R^{3^+2^-}}\cdots. 
$$ 
Note that the operators $\{R^{n^+n^-}\}$ (resp. $\{R^{n^+(n-1)^-}\}$) are commuting among themselves here; 
we may similarly define ${\cal H}_m:=\cdots R^{{m}^+0^-}R^{{m+1}^+1^-}\cdots.$ 
Then the discrete dynamics
\begin{equation}
{\cal T}:=Ad({\cal H}^{-1}) : {\cal O}\mapsto {\cal H}^{-1}{\cal O}{\cal H}, \ {\cal O}\in{\cal A}
\label{def:evolT}
\end{equation}
on $\cal A$ 
is well-defined since any ${\cal O}\in{\cal A}$ is locally supported 
(i.e. of the form $\cdots 1 \otimes a \otimes 1 \cdots$). 
This dynamics 
is explicitly described in terms of 
matrix elements of (finite products of) 
local Lax matrices as we will see shortly (Theorem \ref{thmA}). 

Consider successive products of the local Lax matrices and express them as e.g. 
\begin{eqnarray*}
{\mathcal L_z(1^-1^+2^-2^+)}
&:=&
L_z^-(\Delta^{1-})^{-1} L_z^+(\Delta^{1+}) L_z^-(\Delta^{2-})^{-1} L_z^+(\Delta^{2+})
\end{eqnarray*}
(four points case) and so on. 
If we exploit some more graphics, 
the dynamics ${\cal T}$ applied to 
(the matrix elements of) ${\mathcal L}$ 
can be 
depicted as follows (Fig. 3):
\begin{center}
{\sl Figure 3. The time evolution $\cal T$ (\ref{def:evolT}). }
\end{center}
%

\vskip 0.5cm
\raisebox{4cm}{${\cal L}{
 ({\cdots}1^-1^+2^-2^+{\cdots})}$}
\raisebox{4cm}{=}
\raisebox{3cm}{\includegraphics[height=1.8cm]{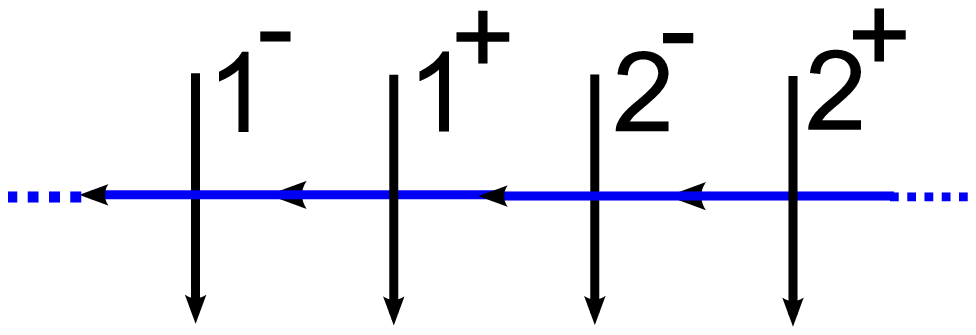}}

\vskip -2cm
\qquad\qquad 
$\mapsto {\cal T(L)}=$
\raisebox{-2.2cm}
{\includegraphics[height=5cm]{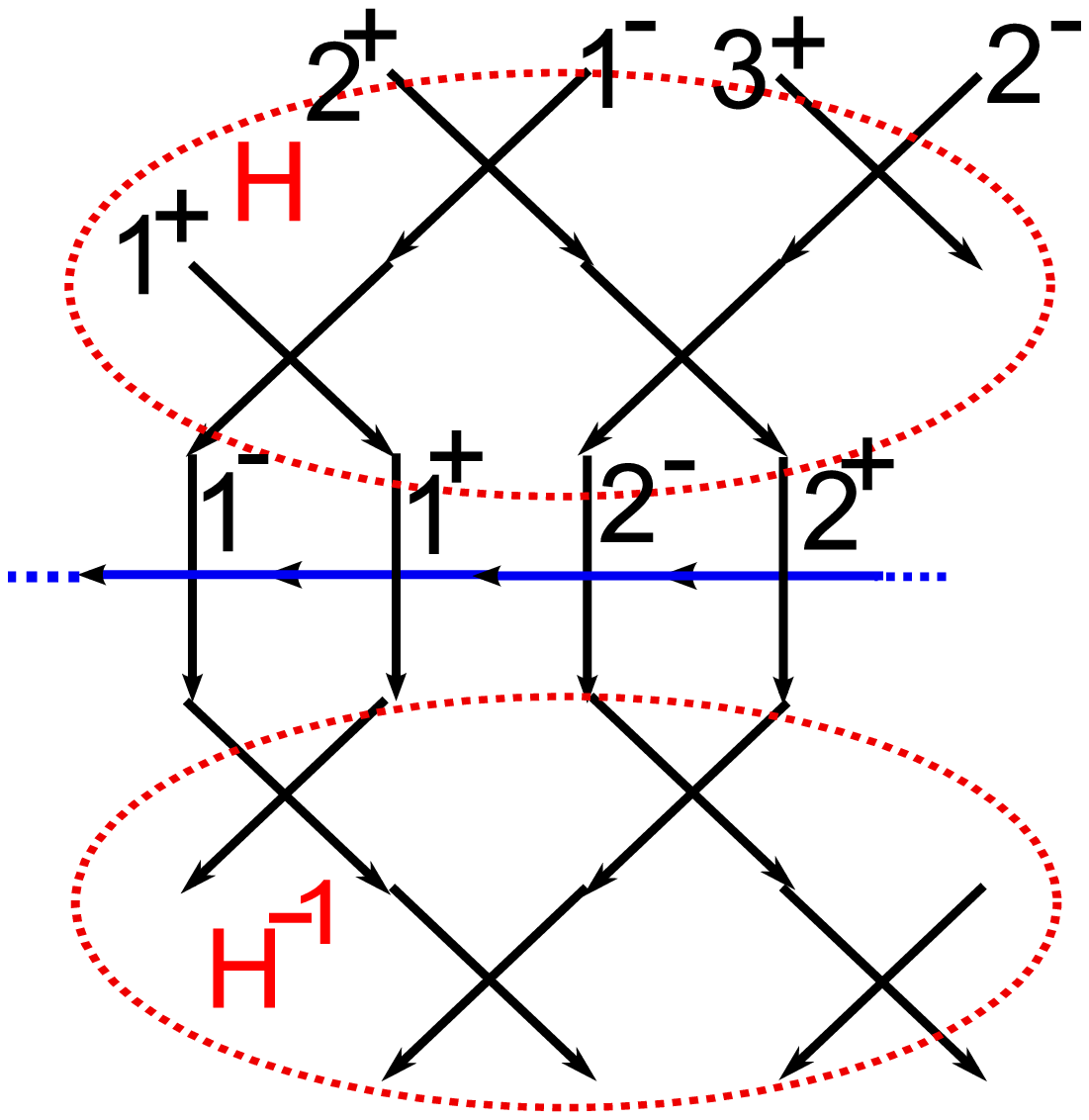}}
\raisebox{3cm}
{\begin{minipage}[t]{5cm}
\quad\quad
\begin{align*}
&{\cal H}_1:=\cdots{R^{2^+1^-}R^{3^+2^-}}\cdots
\\
&{\cal H}_0:=\cdots{R^{1^+1^-}R^{2^+2^-}}\cdots
\\
\quad
\\
&{\cal L}
\\
\quad
\\
&
{\cal H}_0^{-1}
\\
&
{\cal H}_1^{-1}
\end{align*}
\end{minipage}
}

\medskip


To describe the dynamics more explicitly, put 
$\Delta(m^+n^-):=\Delta^{m+}\otimes\Delta^{n-}$ and 
\begin{equation}
w_i(m^+n^-):=
\begin{cases}
\quad (E_{i}^+)^{m+}{\otimes}(E_{i}^-)^{n-} \qquad & (m\equiv n\ {\rm mod}\ 2), 
\\
(E_i^+k_i)^{m+}{\otimes}(k_i^{-1}E_i^-)^{n-} \;\;& \text{(otherwise)}. 
\end{cases}
\label{def:w_i(ab)}
\end{equation}
%
Then 
$
{\cal T}
$
is locally determined by neibouring $w$'s and $\Delta$'s as follows. 
\begin{thm}We have 
\begin{align}
{\cal T}({w_0(1^+1^-)})
&
=
\frac{w_0(1^+0^-){-}q\Delta(1^+0^-)}{w_0(1^+0^-){-}q}
{\cdot}
\frac{w_0(2^+1^-){-}q\Delta(2^+1^-)}{w_0(2^+1^-){-}q}
w_0(2^+0^-)^{-1},
\label{eq:nonautomFV1}
\\
{\cal T}^{-1}({w_0(2^+1^-)})
&=
\frac{w_0(1^+1^-){-}q\Delta(1^+1^-)}{w_0(1^+1^-){-}q}
{\cdot}
\frac{w_0(2^+2^-){-}q\Delta(2^+2^-)}{w_0(2^+2^-){-}q}w_0(1^+2^-)^{-1}.
\label{eq:nonautomFV2}
\end{align}
\label{thmA}
\end{thm}
As for $w_0(m^+, m^-)$ or $w_0((m+1)^+, m^-),$ we read 
$0^\pm, 1^\pm, 2^\pm$ above as $(m-1)^\pm, m^\pm, (m+1)^\pm.$ 
Since (recall (\ref{eq:Deltamoves}))
\begin{equation}
{\cal T}(\Delta^{m\pm})=q^{\mp c^{m\pm}}\Delta^{m\pm},  
\label{eq:evolDelta}
\end{equation}
$
{\cal T}(w_1(m^+,m^-)), {\cal T}(w_1(m^+,{m-1}^-)) 
$
are determined by these formulae. 

\bigskip

{\bf Remark}\ 
From (\ref{eq:localexchangedynamics}), 
we see that ${\cal T}$ induces the exchange dynamics on the lattice :
\begin{align}
{\cal H}_0^{-1}{\cal L}(1^-1^+2^-2^+\cdots n^-n^+){\cal H}_0
&={\cal L}(1^+1^-2^+2^-\cdots n^+n^-)
\label{eq:exchangedynamics1}
\\
{\cal H}_1^{-1}{\cal L}(1^-2^+2^-3^+\cdots n^-{n+1}^+){\cal H}_1
&={\cal L}(2^+1^-3^+2^-\cdots (n+1)^+n^-)
\label{eq:exchangedynamics2}
\end{align}
Unfortunately, the result 
$
{\cal H}_1^{-1}{\cal L}(1^+1^-\cdots n^+n^-)
$
or
$
{\cal H}_0^{-1}{\cal L}(2^+1^-\cdots (n+1)^+n^-) 
$
are not simple enough so that ${\cal T}(\cal L)$ can be said as ``exchange dynamics'' if we 
take the Lax matrix $\cal L$ as a representative of the conjugacy class of
$
{\cal L}\sim L^{-1}{\cal L}L, 
$
together with the following periodic condition. 

\bigskip
Now, let us assume 
 $c^{m\pm}=c^{(m+2)\pm}$
for $\pm=+, -$ and $m=0, 1$ in what follows.
Then we have the trivial $U_q^\pm$-isomorphisms 
$
\iota^{m\pm}: V^{m\pm}\stackrel{\sim}{\rightarrow} V^{(m+2)\pm} \ (1\mapsto 1) ; 
{\rm End}V^{m\pm}\stackrel{\sim}{\rightarrow} {\rm End}V^{m\pm} 
$ 
and therefore 
$$
{\cal S}:=\otimes\iota^{m\pm} : {\cal V}\stackrel{\sim}{\rightarrow}  {\cal V}; \quad
{\cal A}\stackrel{\sim}{\rightarrow}  {\cal A}.
$$
The isomorphism 
${\cal S}$ is nothing but the dilation in the space direction. 
It is obvious that ${\cal H}_0$ and ${\cal H}_1$ do not change with respect to this dilation and hence
\begin{lem}(Periodic reduction) We have
$$
[{\cal H}_0, {\cal S}]=[{\cal H}_1, {\cal S}]=0
$$
so that the dynamics ${\cal T}$ descends to the quotient
$$
\overline{\cal A}:={\cal A}/{\rm Im}({\cal S}-1).
$$
\end{lem}
That is, under the assumption $c_m^{\pm} = c_{m+2}^\pm$, 
the dynamics ${\cal T}$ preserves the conditions
\begin{align}
& w_i(m^+n^-)=
w_i((m+2)^+n^-)=w_i(m^+(n+2)^-), \\
&
\Delta(m^+n^-)=\Delta((m+2)^+n^-)=\Delta(m^+(n+2)^-).
\label{eq:periodicred}
\end{align}

Under this periodic reduction, 
comparison of the obtained formulae 
(\ref{eq:nonautomFV1})
(\ref{eq:nonautomFV2})
with the ones (\ref{qP6^}) in the Appendix via the Weyl group approach, 
we can identify the resulting system with the quantum discrete Painlev\'e system $\widehat{qP_{VI}}$ 
with the following identification of the parameters and the dynamical variables. 
Write 
$w_{i}^{mn}:=w_i(m^+n^-)$ and $\Delta^{mn}:=\Delta(m^+n^-)$for short. 
It turns out that we should identify as follows,
\begin{equation}
-F=
\sqrt[4]{\frac{w_0^{11}w_0^{22}}{w_1^{22}w_1^{11}}}, \
-G=
\sqrt[4]{\frac{w_0^{10}w_0^{21}}{w_1^{21}w_1^{10}}}
\label{eq:identifyF&G}
\end{equation}
and 
$$
a_0^4={\frac{w_0^{21}}{w_0^{10}}}, \
a_1^4={\frac{w_1^{21}}{w_1^{10}}}, \ 
a_2^4=\frac{1}{\Delta^{21}}, \
a_3^4={\Delta^{11}}, \ 
a_4^4={\frac{w_1^{22}}{w_1^{11}}}, \
a_5^4={\frac{w_0^{22}}{w_0^{11}}}.
$$
It is easy to see that 
$a_i$ for $i=0, 1, 4, 5$ and $p=a_0a_1a_2^2a_3^2a_4a_5$ are
central elements among the algebra of observables, 
so that they are constants with respect to our dynamics ${\cal T}$. 
Moreover, 
$$
t:=\Delta(0^+0^-)\Delta(1^+1^-)
$$
satisfies  
${\cal T}(t)=q^{2c}t$, where 
 $c=c^{0-}+c^{1-}-c^{0+}-c^{1+}$ 
(cf. (\ref{eq:evolDelta})), 
meaning that $t$ can be regarded as the time parameter of the dynamics.

\medskip
\begin{thm}
The quantum Painlev\'{e} VI system $\widehat{qP_{VI}}$ (\ref{qP6^}) is reproduced by the above construction:
\begin{equation*}
\begin{cases}
\quad
{\cal T}(F)
&=
\displaystyle
\frac{p^2}{q^2t^2}
\cdot
\frac{G+p^{-1}a_1^{-2}t}{G+p
a_0^{-2}t^{-1}}
\cdot\frac{G+p^{-1}a_1^2t}{G+p
a_0^{2}t^{-1}}F^{-1},
\\
{\cal T}^{-1}(G)
&=
\displaystyle
\frac{1}{q^2t^2}
\cdot
\frac{F+a_4^{2}\,t}{F+a_5^{2}\,t^{-1}}\cdot\frac{F+a_4^{-2}\,t}{F+a_5^{-2}\,t^{-1}}G^{-1}.
\end{cases}
\end{equation*}
\end{thm}

{\bf Remark. }
In (\ref{eq:identifyF&G}) we should employ the fourth root so as to getting the same formula as $\widehat{qP_{IV}}$. 
In fact we can find the $W(D_5^{(1)})$ action without these fourth root and allows us to recover 
$\widehat{qP_{IV}}$ as in the manner in the Appendix.
As in the Faddeev-Volkov system, 
$w_0(1^+1^-), w_0(1^+2^-)$
(or $F$ and $G$) 
together with $a_i^4$ ($i=0, \cdots, 5$) 
generates the diagonal-gauge invariants:
$$\langle w_0(1^+1^-), w_0(1^+2^-), a_0^4, \cdots, a_5^4\rangle
=\langle \{{\cal L}(1^-1^+2^-2^+)_{ij}\ | \ i, j=1, 2\}\rangle^{Ad H}
$$
where 
$H=\left\{
\left[\begin{array}{cc}a & 0\\0 & b\end{array}\right]|a, b\in {\bf C}; a, b \neq 0 
\right\}$. 

\medskip
We also remark that the quantized lattice Liouville equation \cite{FKV} \cite{Kashaev} appears 
as a limit of our equation. 
That is, if we assume 
${\cal T}^N \Delta(1^+0^-), {\cal T}^N \Delta(2^+1^-), {\cal T}^N \Delta(1^+1^-), {\cal T}^N \Delta(2^+2^-)\rightarrow 0$ 
as $N\rightarrow \infty$ (i.e. ${\rm Re} (c^{m+}-c^{n-})>0$ for all $m, n$), then 
from (\ref{eq:nonautomFV1}), (\ref{eq:nonautomFV2}) we respectively have
\begin{align*}
{\cal T}({w_0(1^+1^-)})
w_0(2^+0^-)
&
=
\frac{w_0(1^+0^-)
}{w_0(1^+0^-){-}q}
{\cdot}
\frac{w_0(2^+1^-)
}{w_0(2^+1^-){-}q},
\\
{\cal T}^{-1}({w(2^+1^-)})
w_0(1^+2^-)
&=
\frac{w_0(1^+1^-)
}{w_0(1^+1^-){-}q}
{\cdot}
\frac{w_0(2^+2^-)
}{w_0(2^+2^-){-}q}
\end{align*}
or
\begin{equation}
\begin{cases}
&
w_0(2^+0^-)^{-1}
{\cal T}({w_0(1^+1^-)}^{-1})
=
\left(1-qw_0(1^+0^-)^{-1}\right)
\left(1-qw_0(2^+1^-)^{-1}\right),
\\
& w_0(1^+2^-)^{-1}
{\cal T}^{-1}({w(2^+1^-)}^{-1})
=
\left(1-qw_0(1^+1^-)^{-1}\right)
\left(1-qw_0(2^+2^-)^{-1}\right).
\end{cases}
\label{Liouville}
\end{equation}

%

\section{Quantized 
discrete Schlesinger equation}

\begin{thm}
Let ${\cal B}(z):={\cal H}_0{\cal L}_z(2^-){\cal H}_1 $. We have  
\begin{equation}
{\cal L}(1^-1^+2^-2^+)
{\cal B}(z)
=
{\cal B}(z)
{\cal L}(1^+2^-2^+1^-).
\label{eq:thmC}
\end{equation}
\label{thm:C}
\end{thm}
Let us write
$
{\cal L}(1^-1^+2^-2^+)=
q^{D}L(1^-1^+2^-2^+), \ 
$ 
where 
$q^{D}$ stands for the 
difference operator part, 
$$
D
=\sum_{i=1,2;\pm}\mp c^{i\pm} d =cd.
$$
Then the above relation (\ref{eq:thmC}) is equivalent to (cf. (\ref{eq:homogeneity}))
\begin{equation}
L(1^-1^+2^-2^+){\cal B}(z)={\cal B}(zq^{-c}){\cal T}(L(1^-1^+2^-2^+)), 
\label{eq:thmC'}
\end{equation}
which can be recognized as the 
quantization of (\ref{eq:qSchlesinger}). 

{\bf Proof} 
\begin{align*}
{\rm LHS}&={\cal L}(1^-1^+2^-2^+){\cal H}_0 L(2^-){\cal H}_1
\\
&={\cal H}_0{\cal L}(1^+1^-2^+2^-) L(2^-){\cal H}_1
\\
&={\cal H}_0L(2^-){\cal L}(2^-1^+1^-2^+){\cal H}_1
\\
&={\cal H}_0L(2^-){\cal H}_1
{\cal L}(1^+2^-2^+1^-)={\rm RHS}.\qquad\qquad\framebox{}
\end{align*}
The above proof uses 
the Yang-Baxter equation under the periodicity condition, 
which can be depicted as Fig. 4. 
\begin{center}
{\sl Figure 4. The discrete Schlesinger equation (\ref{eq:thmC}). }
\\
\bigskip
\includegraphics[height=4cm]{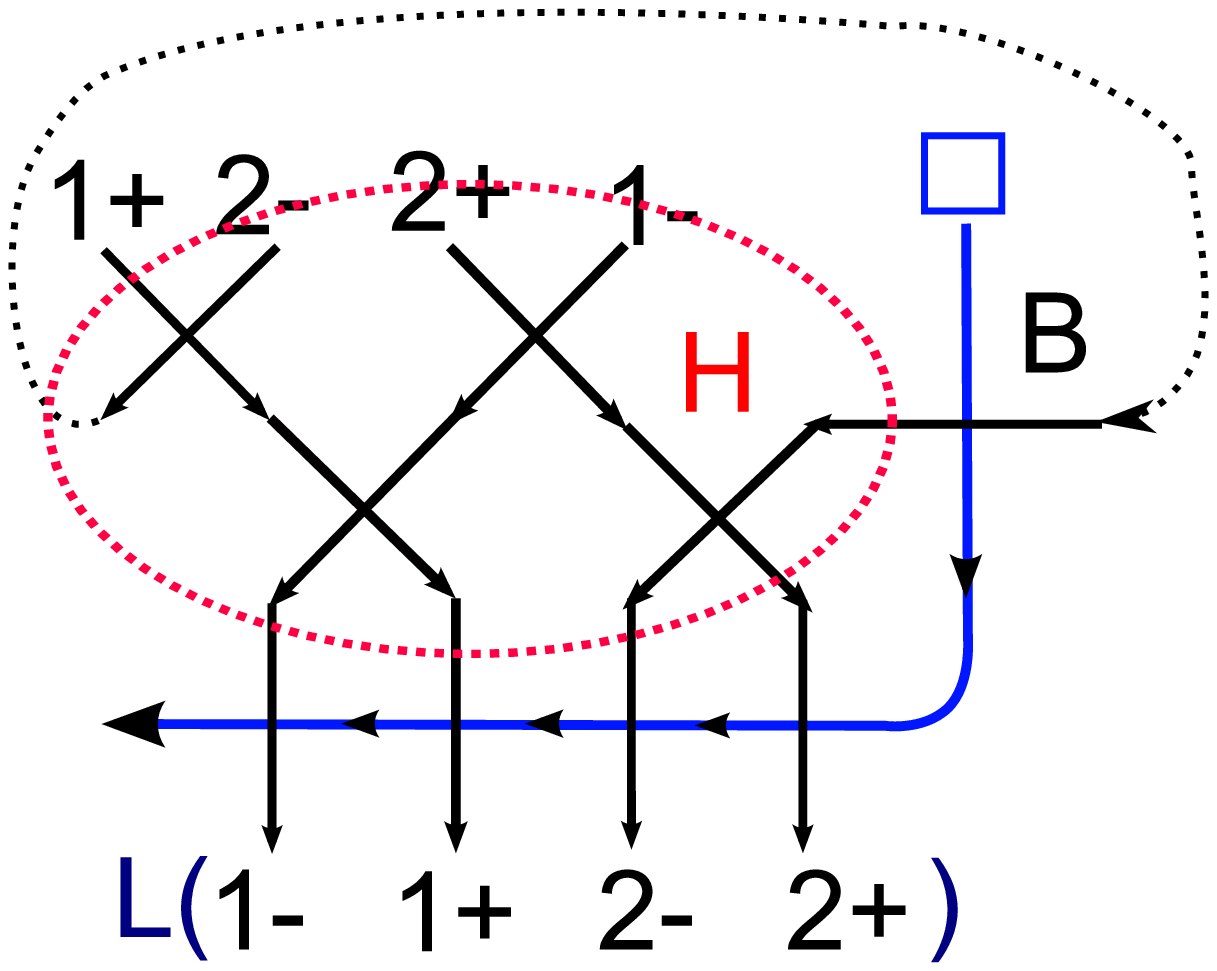}
\raisebox{1.5cm}{$=$}
\includegraphics[height=4cm]{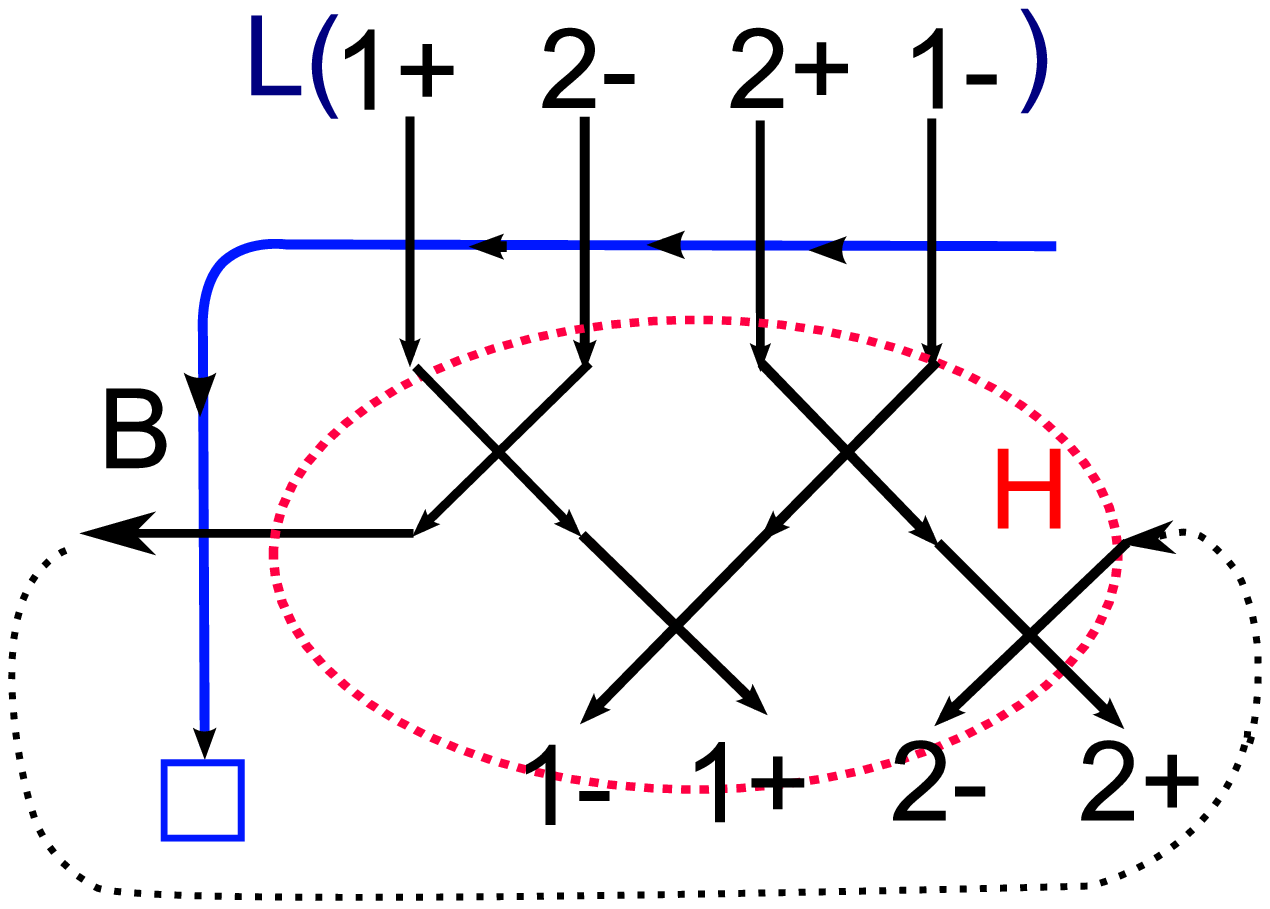}
\end{center}

\bigskip
In general, without the periodicity, we have 
(cf. (\ref{eq:exchangedynamics1})(\ref{eq:exchangedynamics2}))
\begin{equation}
{\cal L}(1^-1^+\cdots n^- n^+){\cal B}_n(z)={\cal B}_0(z){\cal L}(1^+0^-2^+1^-\cdots n^+(n-1)^-)
\label{eq:LB=BL'ingeneral}
\end{equation}
on ${\cal V}$, where 
$
{\cal B}_k(z):={\cal H}_0 L(k^-){\cal H}_1 
$
for $k=0, n$. 
If we assume the $n$- periodicity 
$c^{m\pm}=c^{(m+n)\pm}, V^{m\pm}\stackrel{\sim}{\rightarrow} V^{(m+n)\pm}$, then 
$
{\cal L}(1^+0^-2^+1^-\cdots n^+(n-1)^-)
=
{\cal L}(1^+n^-2^+1^-\cdots n^+(n-1)^-)
$ 
in the right-hand side is conjugate to
${\cal T}({\cal L}(1^-1^+\cdots n^- n^+))$, 
so that the equation 
(\ref{eq:LB=BL'ingeneral})
((\ref{eq:thmC}) or (\ref{eq:thmC'}) in $n=2$)
can be considered as the compatibility condition of the 
linear problem 
\begin{equation}
{\cal Y}={\cal L}{\cal Y}, \ {\cal L}={\cal L}(1^-1^+\cdots n^- n^+), 
\label{eq:linprob}
\end{equation}
or, equivalently (write ${\cal L}=q^{D}L(z)$ as in (\ref{eq:thmC'}) )
\begin{equation}
q^{-D}{\cal Y}(z)={\cal Y}(q^{-c}z)=L(z){\cal Y}(z), 
\label{eq:linprob^}
\end{equation}
(where 
$
{\cal Y}
$
should be regarded as 
$\overline
{\cal V}
\otimes {\bf C}^2$- valued) 
and the time evotution ${\cal T}.$

Thus we have succeeded in quantizing the isomonodromy problem 
or the Lax form for the quantum discrete Painlev\'e system $\widehat{qP_{VI}}$ explicitly.

It is interesting to note that two fundamental solutions of (\ref{eq:linprob^}) can be
at least formally obtained as 
$$
{\cal Y}_\infty(z):=L(z)L(zq^c)L(zq^{2c})\cdots
$$
and 
$$
{\cal Y}_0(z):=L(zq^{-c})^{-1}L(zq^{-2c})^{-1}\cdots
$$
and then the quantization of (\ref{eq:connectionmatrix}) can be simply written as
$$
{\cal M}(z):=
{\cal Y}_0^{-1}{\cal Y}_\infty
=\cdots L(zq^{-2c})L(zq^{-c})L(z)L(zq^{c})L(zq^{2c})\cdots.
$$

\section{
The Weyl group action}

Further comparison with the formula in the Appendix give us the 
formula for the $\tilde{W}(D_5^{(1)})=\langle W, \sigma\rangle$- action. 
We have
$$
\sigma: w_0^{11}\mapsto q^{-2}w_1^{11}, \ w_0^{21}\mapsto q^{-2}w_1^{21}, 
$$
and
$$
s_0: w_0^{11}\mapsto w_0^{11}, \ w_0^{21}\mapsto w_0^{10}, \qquad\qquad\qquad\qquad\qquad\qquad
s_5: w_0^{11}\mapsto w_0^{22}, \ w_0^{21}\mapsto w_0^{21}, 
$$
$$
s_2: w_0^{11}\mapsto w_0^{11}\frac{w_0^{21}-1}{w_0^{21}-\Delta^{21}}, \ w_0^{21}\mapsto \frac{1}{w_1^{21}}, 
\ 
s_3: w_0^{11}\mapsto \frac{1}{w_1^{11}}, \ w_0^{21}\mapsto \frac{w_0^{11}-1}{w_0^{11}-\Delta^{11}}w_0^{10}, 
$$
$$
s_1: w_0^{11}\mapsto w_0^{11}, \ w_0^{21}\mapsto w_0^{21}, \qquad\qquad\qquad\qquad\qquad\qquad
s_4: w_0^{11}\mapsto w_0^{11}, \ w_0^{21}\mapsto w_0^{21}.
$$

\bigskip
There should be a Lax matrix point of view elucidation of these symmetry : 
it is quite plausible that these symmetry arise from the choice of 
multiplicative decompositions of our Lax matrix (cf. \cite{Borodin}), 
and is related to the tesseretion of the projective line with four points 
so that the time evolution ${\cal T}$ can be regarded as the Dehn twist 
(\cite{Kashaev98}). 
We would like to report this point seperately 
in a near future. 

\bigskip

{\bf Acknowledgement.} 
The author is supported by the Grant-in-Aid for Scientific Research (Kakenhi) C, 23540004. 
Part of this paper is based on a lecture dilivered at University of Tokyo, 2007, and the author is 
greateful for Professor Jimbo Michio and Professor Jun-ichi Shiraishi for their hospitality.
Thanks are also due to 
 Professor Gen Kuroki, Professor Yasuhiko Yamada, 
and (especially) 
Professor Akihiro Tsuchiya 
for their kind interest.

\section
{Appendix. Review of the Weyl group approach to {\bf $\widehat{qP_{VI}}$}
\cite{H}}

Here we recite the results of \cite{H}. 
We use
the 
$D_5^{(1)}$ root system. 
Let $\{e_j\}_{1\leq j\leq 5}$ be the orthonormal basis $\subset {\bm R}^6={\bm R}^5 \perp{\bm R}\delta$, 
where $\delta$ is identified with the null root, then the symple roots will be realized as follows. 

\bigskip
\quad
\begin{small}
$
\begin{array}{cccccccc}
\alpha_0 &  &   &   &   &  & \alpha_5\\
     &\backslash&   &   &   &\slash&  \\
     &  &\alpha_2 & - &\alpha_3 &  &  \\
  &\slash&   &   &   &\backslash&  \\
\alpha_1 &  &   &   &   &  & \alpha_4
\end{array}
=
\begin{array}{rccccl}
{\small\delta-e_1-e_2} &   &   &   &  {\small e_4+e_5}\\
                     \backslash&   &   &   &\slash\qquad\\
          &{\small e_2-e_3}& - & {\small e_3-e_4} &    \\
                     \slash&   &   &   &\backslash\qquad \\
     {\small e_1-e_2}  &   &   &   &  {\small e_4-e_5}
\end{array}
$
\end{small}
\bigskip


Let $q=e^\hbar\in {\bf C}^\times, |q|<1$. 
Let $a_j:=e^{\hbar\alpha_j}, {p}:=e^{\hbar\delta}=a_0a_1a_2^2a_3^2a_4a_5$ be elements of the 
group algebra of the $D_5^{(1)}$ root lattice.   
We have the 
$
 W=W(D_5^{(1)})$ action given by

$$
s_i(a_j)=a_i^{-C_{ij}}a_j
\quad (s_i(p)=p, \  \forall i).
$$
%
We also need diagram automorphisms. They are
$$
\tau
: a_j\longleftrightarrow a_{5-j}^{-1} (j=0,\cdots,5),
$$
and
$$
\sigma: a_0\leftrightarrow a_1^{-1}, 
\ 
a_4\leftrightarrow a_5^{-1}, a_j\mapsto a_j^{-1} (j=2,3). 
$$
We have defined the action of the extended affine Weyl group 
 $\tilde{W}=\langle W, \tau, \sigma 
\rangle$ on the group algebra of the root lattice $Q$.
Let further 
 ${\bm K}:=\bm C(a_0, \cdots, a_5)\langle F,  G\rangle$
where $F$ and $G$ are noncommutative letters; we let {$FG=q^2GF$} later.

\bigskip
{\bf Theorem}\ (1) 
We have 
 $\tilde W(D_5^{(1)})$-action

$$
 \langle W, \sigma
\rangle\stackrel{hom}{\rightarrow} {\rm Aut}(\bm K)
$$
given by
$$
\sigma: F \leftrightarrow q^{-2}F^{-1}, G \leftrightarrow q^{-2}G^{-1},
$$

$$
{s_2}(F):=  F \frac{a_0a_1^{-1}G+a_2^2}{a_0a_1^{-1}a_2^2G+1}, \quad
s_j(F):= F \quad (j\neq 2)
$$
$$
{s_3}(G):=   \frac{a_3^2a_4a_5^{-1}F+1}{a_4a_5^{-1}F+a_3^2}  G, \quad
s_j(G):= G \quad (j\neq 3)
$$

(2)
If {$FG=q^2GF$}, this 
action is {Hamiltonian}:
namely we have
$
{\Sigma, S_j}$ such that 
$$
{\sigma_{01}\sigma_{45}(}\phi{)=\Sigma} \phi {\Sigma^{-1}},
\quad
{s_j(}\phi{)=S_j}\phi {S_j^{-1}} 
 (j=0, \cdots, 5)
$$
for any $\phi \in \bm K$. 
Recall 
$
{(x)_\infty}=(x,q)_\infty 
= \prod_{m=0}^\infty (1-xq^m)$ and put 
$
{\Psi(a,x)}
:=
\frac
{(qx)_\infty (x^{-1})_\infty}
{(aqx)_\infty (ax^{-1})_\infty}.
$
Then $$
{\Sigma}
:=
(FG)_\infty 
(qG^{-1}F^{-1})_\infty 
(G^{-1}F)_\infty
(qF^{-1}G)_\infty 
(F)_\infty^2
(qF^{-1})_\infty^2
(G)_\infty^2
(qG^{-1})_\infty^2,
$$
and
$$
{S_2}:=\Psi(a_2, a_0a_1^{-1}G){\rm e}^{\frac{\pi i}{2}\alpha_2\partial_2}, 
{S_3}:=\Psi(a_3, a_5a_4^{-1}G){\rm e}^{\frac{\pi i}{2}\alpha_3\partial_3},
{S_j}
:=
{\rm e}^{\frac{\pi i}{2}\alpha_j\partial_j} (j\neq 2,3),
$$
where the derication $\partial_j$ is defined by 
$
{\partial_j}(\alpha_k):=C_{j,k}
$
(the Cartan matrix).
\hfill\framebox{}



\medskip



Now consider the lattice element 
$T_3:=s_2s_1s_0s_2s_3s_4s_5s_3\sigma_{01}\sigma_{45}
: 
e_j\mapsto e_j-\delta_{j,3}\delta, \; \delta\mapsto\delta.
$
In $q=1$ case, $T_3$ 
reproduces the {$qP_{VI}$} of Jimbo-Sakai 
\cite{TM}. 
Put $t=q^{2e_3}=a_3^2a_4a_5$. 

\bigskip
{{\bf Theorem/Definition.}}
The $T_3$ action is given as follows ({\sl the quantized difference Painlev\'e VI system} $\widehat{qP_{VI}}$), 
which commutes with 
$
W(D_4^{(1)})
\simeq 
\langle s_0, s_1, s_2s_3s_2, s_4, s_5\rangle. 
$
$$
T_3(a_0, a_1, t, a_4, a_5)=
(a_0, a_1, t/p, a_4, a_5),
$$
\begin{equation}
{\bf \widehat{qP_{VI}}} :
\begin{cases}
\ T_3(F)
&=
%
\displaystyle
\frac{p^2}{q^2t^2}
\cdot
\frac{G+p^{-1}a_1^{-2}t}{G+p
a_0^{-2}t^{-1}}
\cdot\frac
{G+p^{-1}a_1^2t}{G+p
a_0^{2}t^{-1}}F^{-1}
, 
\\ 
\ T_3^{-1}(G)
&=
\displaystyle
\frac{1}{q^2t^2}
\cdot
\frac{F+a_4^{2}\,t}{F+a_5^{2}\,t^{-1}}
\cdot\frac{F+a_4^{-2}\,t}{F+a_5^{-2}\,t^{-1}}
G^{-1}.
\end{cases}
\label{qP6^}
\end{equation}
If 
$FG=q^2GF$, then by construction, $\widehat{qP_{VI}}(= T_3$ action) has the Hamiltonian, 
$$
{T_3=Ad({\cal H})}, \ {{\cal H}:=S_2S_1S_0S_2S_3S_4S_5S_3\Sigma}.
$$



\newpage

\end{document}